# The Lie algebra *split*$g_2$ with *Mathematica* using Zorn's matrices.

**Pablo Alberca Bjerregaard**

*Department of Applied Mathematics*
*Escuela Técnica Superior de Ingeniería Industrial*
*Campus de El Ejido*
*University of Málaga*
*Málaga –Spain*
*E–Mail:pgalberca@uma.es*

**Cándido Martín González**

*Department of Algebra, Geometry and Topology*
*Facultad de Ciencias*
*University of Málaga*
*Campus de Teatinos*
*Málaga –Spain*
*E–Mail:candido@apncs.cie.uma.es*

*2001.*

**Abstract.** *We will obtain in this paper a generic expression of any element in* split $g_2$ *,the Lie algebra of the derivations of the split octonions* $\mathbf{O}_s$ *over an arbitrary field. For this purpose, we will use the Zorn's matrices. We will also compute the multiplication table of this Lie algebra.*

## ■ Previous commands

First of all, we define some commands that will be usefull in what follows. The first ones are **vec**, vectorial product, **esc**, scalar product and **prod**, the Zorn's matrices product:

```
vec[{x_, y_, z_}, {u_, v_, w_}] =
    {y*w - z*v, z*u - x*w, x*v - y*u};
esc[{x_, y_, z_}, {u_, v_, w_}] =    x*u + y*v + z*w;
prod[( α_  x_ ), ( γ_  z_ )] :=
     ( y_  β_ )  ( t_  δ_ )
    (     α*γ + esc[x, t]        α*z + δ*x - vec[y, t] )
    ( γ*y + β*t + vec[x, z]        β*δ + esc[y, z]     )
```

The function **Recon** writes the matrix of any derivation as a linear combination of the future basis elements $x_j$ :

```
Recon[m_] := Module[{},
    λ1 = m[[1, 3]]; λ2 = m[[1, 4]]; λ3 = m[[1, 5]];
    λ4 = m[[3, 3]]; λ5 = m[[3, 4]]; λ6 = m[[3, 5]];
    λ7 = m[[4, 3]]; λ8 = m[[4, 4]]; λ9 = m[[4, 5]];
    λ10 = m[[5, 3]]; λ11 = m[[5, 4]]; λ12 = m[[1, 6]];
    λ13 = m[[1, 7]]; λ14 = m[[1, 8]]; Sum[λi * (HoldForm[x])i , {i, 14}]]
```

Next we define the commands **c**, Lie product, and **δ** the function Krönecker delta:

```
c[x_, y_] := x.y - y.x;
δi_,j_ := If[i == j, 1, 0]
```

## ■ Generic Matrix

We define now the set {$\Theta$, $\Upsilon$, A, B, $C_1$, $C_2$, $C_3$, $D_1$, $D_2$, $D_3$} where the last 8 Zorn's matrices constitute a basis of $\mathbf{O}_s$ .

$$\theta = \begin{pmatrix} 0 & \{0,0,0\} \\ \{0,0,0\} & 0 \end{pmatrix}; \Upsilon = \begin{pmatrix} 1 & \{0,0,0\} \\ \{0,0,0\} & 1 \end{pmatrix};$$
$$A = \begin{pmatrix} 1 & \{0,0,0\} \\ \{0,0,0\} & 0 \end{pmatrix}; B = \begin{pmatrix} 0 & \{0,0,0\} \\ \{0,0,0\} & 1 \end{pmatrix};$$
$$C_1 = \begin{pmatrix} 0 & \{1,0,0\} \\ \{0,0,0\} & 0 \end{pmatrix}; C_2 = \begin{pmatrix} 0 & \{0,1,0\} \\ \{0,0,0\} & 0 \end{pmatrix};$$
$$C_3 = \begin{pmatrix} 0 & \{0,0,1\} \\ \{0,0,0\} & 0 \end{pmatrix}; D_1 = \begin{pmatrix} 0 & \{0,0,0\} \\ \{1,0,0\} & 0 \end{pmatrix};$$
$$D_2 = \begin{pmatrix} 0 & \{0,0,0\} \\ \{0,1,0\} & 0 \end{pmatrix}; D_3 = \begin{pmatrix} 0 & \{0,0,0\} \\ \{0,0,1\} & 0 \end{pmatrix};$$

If $d : \mathbf{O}_s \longrightarrow \mathbf{O}_s$ is a derivation, we are going to determine its associtated matrix in the previous basis. We write the images of the basis elements as generic Zorn's matrices using parameters:

$$dA = \begin{pmatrix} a_1 & \{u_{11}, u_{12}, u_{13}\} \\ \{v_{11}, v_{12}, v_{13}\} & b_1 \end{pmatrix}; dB = \begin{pmatrix} a_2 & \{u_{21}, u_{22}, u_{23}\} \\ \{v_{21}, v_{22}, v_{23}\} & b_2 \end{pmatrix};$$
$$dC_1 = \begin{pmatrix} a_3 & \{u_{31}, u_{32}, u_{33}\} \\ \{v_{31}, v_{32}, v_{33}\} & b_3 \end{pmatrix}; dC_2 = \begin{pmatrix} a_4 & \{u_{41}, u_{42}, u_{43}\} \\ \{v_{41}, v_{42}, v_{43}\} & b_4 \end{pmatrix};$$
$$dC_3 = \begin{pmatrix} a_5 & \{u_{51}, u_{52}, u_{53}\} \\ \{v_{51}, v_{52}, v_{53}\} & b_5 \end{pmatrix}; dD_1 = \begin{pmatrix} a_6 & \{u_{61}, u_{62}, u_{63}\} \\ \{v_{61}, v_{62}, v_{63}\} & b_6 \end{pmatrix};$$
$$dD_2 = \begin{pmatrix} a_7 & \{u_{71}, u_{72}, u_{73}\} \\ \{v_{71}, v_{72}, v_{73}\} & b_7 \end{pmatrix}; dD_3 = \begin{pmatrix} a_8 & \{u_{81}, u_{82}, u_{83}\} \\ \{v_{81}, v_{82}, v_{83}\} & b_8 \end{pmatrix};$$

We will determine these parameters using the derivation condition and the relations between them. Firstly, as

```
A + B == Υ
```

True

we apply the derivation $d$ to obtain that dB = −dA , and then

$$dA = \{\{a_1, \{u_{11}, u_{12}, u_{13}\}\}, \{\{v_{11}, v_{12}, v_{13}\}, b_1\}\};$$
$$dB = \{\{-a_1, \{-u_{11}, -u_{12}, -u_{13}\}\}, \{\{-v_{11}, -v_{12}, -v_{13}\}, -b_1\}\};$$
$$dC_1 = \{\{a_3, \{u_{31}, u_{32}, u_{33}\}\}, \{\{v_{31}, v_{32}, v_{33}\}, b_3\}\};$$
$$dC_2 = \{\{a_4, \{u_{41}, u_{42}, u_{43}\}\}, \{\{v_{41}, v_{42}, v_{43}\}, b_4\}\};$$
$$dC_3 = \{\{a_5, \{u_{51}, u_{52}, u_{53}\}\}, \{\{v_{51}, v_{52}, v_{53}\}, b_5\}\};$$
$$dD_1 = \{\{a_6, \{u_{61}, u_{62}, u_{63}\}\}, \{\{v_{61}, v_{62}, v_{63}\}, b_6\}\};$$
$$dD_2 = \{\{a_7, \{u_{71}, u_{72}, u_{73}\}\}, \{\{v_{71}, v_{72}, v_{73}\}, b_7\}\};$$
$$dD_3 = \{\{a_8, \{u_{81}, u_{82}, u_{83}\}\}, \{\{v_{81}, v_{82}, v_{83}\}, b_8\}\};$$

As $A = A\,A$

```
A == prod[A, A]
```

True

we have that

```
-dA + prod[dA, A] + prod[A, dA]
```

$$\begin{pmatrix} a_1 & \{0, 0, 0\} \\ \{0, 0, 0\} & -b_1 \end{pmatrix}$$

which must be null. We obtain then $a_1 = b_1 = 0$, and

$$dA = \{\{0, \{u_{11}, u_{12}, u_{13}\}\}, \{\{v_{11}, v_{12}, v_{13}\}, 0\}\};$$
$$dB = \{\{0, \{-u_{11}, -u_{12}, -u_{13}\}\}, \{\{-v_{11}, -v_{12}, -v_{13}\}, 0\}\};$$
$$dC_1 = \{\{a_3, \{u_{31}, u_{32}, u_{33}\}\}, \{\{v_{31}, v_{32}, v_{33}\}, b_3\}\};$$
$$dC_2 = \{\{a_4, \{u_{41}, u_{42}, u_{43}\}\}, \{\{v_{41}, v_{42}, v_{43}\}, b_4\}\};$$
$$dC_3 = \{\{a_5, \{u_{51}, u_{52}, u_{53}\}\}, \{\{v_{51}, v_{52}, v_{53}\}, b_5\}\};$$
$$dD_1 = \{\{a_6, \{u_{61}, u_{62}, u_{63}\}\}, \{\{v_{61}, v_{62}, v_{63}\}, b_6\}\};$$
$$dD_2 = \{\{a_7, \{u_{71}, u_{72}, u_{73}\}\}, \{\{v_{71}, v_{72}, v_{73}\}, b_7\}\};$$
$$dD_3 = \{\{a_8, \{u_{81}, u_{82}, u_{83}\}\}, \{\{v_{81}, v_{82}, v_{83}\}, b_8\}\};$$

From the identities

```
prod[C1, C1] == θ
prod[C2, C2] == θ
prod[C3, C3] == θ
```

True

True

True

we conclude that the following matrices must also be null

```
prod[dC1, C1] + prod[C1, dC1]
prod[dC2, C2] + prod[C2, dC2]
prod[dC3, C3] + prod[C3, dC3]
```

$$\begin{pmatrix} v_{31} & \{a_3 + b_3, 0, 0\} \\ \{0, 0, 0\} & v_{31} \end{pmatrix}$$

$$\begin{pmatrix} v_{42} & \{0, a_4 + b_4, 0\} \\ \{0, 0, 0\} & v_{42} \end{pmatrix}$$

$$\begin{pmatrix} v_{53} & \{0, 0, a_5 + b_5\} \\ \{0, 0, 0\} & v_{53} \end{pmatrix}$$

and then $v_{31} = v_{42} = v_{53} = 0$, $b_3 = -a_3$, $b_4 = -a_4$ y $b_5 = -a_5$, and so we have

```
dA = {{0, {u11, u12, u13}}, {{v11, v12, v13}, 0}};
dB = {{0, {-u11, -u12, -u13}}, {{-v11, -v12, -v13}, 0}};
dC1 = {{a3, {u31, u32, u33}}, {{0, v32, v33}, -a3}};
dC2 = {{a4, {u41, u42, u43}}, {{v41, 0, v43}, -a4}};
dC3 = {{a5, {u51, u52, u53}}, {{v51, v52, 0}, -a5}};
dD1 = {{a6, {u61, u62, u63}}, {{v61, v62, v63}, b6}};
dD2 = {{a7, {u71, u72, u73}}, {{v71, v72, v73}, b7}};
dD3 = {{a8, {u81, u82, u83}}, {{v81, v82, v83}, b8}};
```

Using the same argument, the identities

```
prod[C1, C2] + prod[C2, C1] == θ
prod[C1, C3] + prod[C3, C1] == θ
prod[C2, C3] + prod[C3, C2] == θ
```

True

True

True

justify that the Zorn's matrices

```
prod[dC1, C2] + prod[C1, dC2] + prod[dC2, C1] + prod[C2, dC1]
prod[dC1, C3] + prod[C1, dC3] + prod[dC3, C1] + prod[C3, dC1]
prod[dC2, C3] + prod[C2, dC3] + prod[dC3, C2] + prod[C3, dC2]
```

$$\begin{pmatrix} v_{32}+v_{41} & \{0,0,0\} \\ \{0,0,0\} & v_{32}+v_{41} \end{pmatrix}$$

$$\begin{pmatrix} v_{33}+v_{51} & \{0,0,0\} \\ \{0,0,0\} & v_{33}+v_{51} \end{pmatrix}$$

$$\begin{pmatrix} v_{43}+v_{52} & \{0,0,0\} \\ \{0,0,0\} & v_{43}+v_{52} \end{pmatrix}$$

must be null. We obtain then that $-v_{32}=v_{41}$, $-v_{33}=v_{51}$ and $-v_{43}=v_{52}$. We have then

```
dA = {{0, {u11, u12, u13}}, {{v11, v12, v13}, 0}};
dB = {{0, {-u11, -u12, -u13}}, {{-v11, -v12, -v13}, 0}};
dC1 = {{a3, {u31, u32, u33}}, {{0, v32, v33}, -a3}};
dC2 = {{a4, {u41, u42, u43}}, {{-v32, 0, v43}, -a4}};
dC3 = {{a5, {u51, u52, u53}}, {{-v33, -v43, 0}, -a5}};
dD1 = {{a6, {u61, u62, u63}}, {{v61, v62, v63}, b6}};
dD2 = {{a7, {u71, u72, u73}}, {{v71, v72, v73}, b7}};
dD3 = {{a8, {u81, u82, u83}}, {{v81, v82, v83}, b8}};
```

The identity $C_1\ (C_2\ C_3) = A$

```
A == prod[C1, prod[C2, C3]]
```

True

implies that the following matrix must vanish

```
-dA + prod[dC1, prod[C2, C3]] + prod[C1, prod[dC2, C3]] +
        prod[C1, prod[C2, dC3]]
```

$$\begin{pmatrix} u_{31}+u_{42}+u_{53} & \{v_{43}-u_{11}, -u_{12}-v_{33}, v_{32}-u_{13}\} \\ \{-a_3-v_{11}, -a_4-v_{12}, -a_5-v_{13}\} & 0 \end{pmatrix}$$

We obtain the relations $u_{53}=-u_{31}-u_{42}$, $v_{43}=u_{11}$, $v_{33}=-u_{12}$, $v_{32}=u_{13}$, $a_3=-v_{11}$, $a_4=-v_{12}$ and $a_5=-v_{13}$. Then we have

```
dA = {{0, {u11, u12, u13}}, {{v11, v12, v13}, 0}};
dB = {{0, {-u11, -u12, -u13}}, {{-v11, -v12, -v13}, 0}};
dC1 = {{-v11, {u31, u32, u33}}, {{0, u13, -u12}, v11}};
dC2 = {{-v12, {u41, u42, u43}}, {{-u13, 0, u11}, v12}};
dC3 = {{-v13, {u51, u52, -u31 - u42}}, {{u12, -u11, 0}, v13}};
dD1 = {{a6, {u61, u62, u63}}, {{v61, v62, v63}, b6}};
dD2 = {{a7, {u71, u72, u73}}, {{v71, v72, v73}, b7}};
dD3 = {{a8, {u81, u82, u83}}, {{v81, v82, v83}, b8}};
```

We do not obtain any information from the identity $(C_1\ C_2)\ C_3$

```
B == prod[prod[C1, C2], C3]
```

True

because the matrix

```
-dB + prod[prod[dC1, C2], C3] + prod[prod[C1, dC2], C3] + prod[prod[C1, C2], dC3]
```

$$\begin{pmatrix} 0 & \{0,0,0\} \\ \{0,0,0\} & 0 \end{pmatrix}$$

is directly null. As $C_2\ C_3 = D_1$

```
D1 == prod[C2, C3]
```

True

the matrix

```
-dD1 + prod[dC2, C3] + prod[C2, dC3]
```

$$\begin{pmatrix} -a_6 - u_{11} & \{-u_{61}, v_{13} - u_{62}, -u_{63} - v_{12}\} \\ \{-u_{31} - v_{61}, -u_{41} - v_{62}, -u_{51} - v_{63}\} & u_{11} - b_6 \end{pmatrix}$$

must be null. Then we have $a_6 = -u_{11}$, $b_6 = u_{11}$, $u_{61} = 0$, $u_{62} = v_{13}$, $u_{63} = -v_{12}$, $v_{61} = -u_{31}$, $v_{62} = -u_{41}$ and $v_{63} = -u_{51}$. Then we write

```
dA = {{0, {u11, u12, u13}}, {{v11, v12, v13}, 0}};
dB = {{0, {-u11, -u12, -u13}}, {{-v11, -v12, -v13}, 0}};
dC1 = {{-v11, {u31, u32, u33}}, {{0, u13, -u12}, v11}};
dC2 = {{-v12, {u41, u42, u43}}, {{-u13, 0, u11}, v12}};
dC3 = {{-v13, {u51, u52, -u31 - u42}}, {{u12, -u11, 0}, v13}};
dD1 = {{-u11, {0, v13, -v12}}, {{-u31, -u41, -u51}, u11}};
dD2 = {{a7, {u71, u72, u73}}, {{v71, v72, v73}, b7}};
dD3 = {{a8, {u81, u82, u83}}, {{v81, v82, v83}, b8}};
```

As $D_2 = -C_1\ C_3$

```
D2 == -prod[C1, C3]
```

True

then

```
dD2 + prod[dC1, C3] + prod[C1, dC3]
```

$$\begin{pmatrix} a_7 + u_{12} & \{u_{71} + v_{13}, u_{72}, u_{73} - v_{11}\} \\ \{u_{32} + v_{71}, u_{42} + v_{72}, u_{52} + v_{73}\} & b_7 - u_{12} \end{pmatrix}$$

must be zero. Thus we have $a_7 = -u_{12}$, $b_7 = u_{12}$, $u_{71} = -v_{13}$, $u_{72} = 0$, $u_{73} = v_{11}$, $v_{71} = -u_{32}$, $v_{72} = -u_{42}$ and $v_{73} = -u_{52}$. We have then

```
dA = {{0, {u11, u12, u13}}, {{v11, v12, v13}, 0}};
dB = {{0, {-u11, -u12, -u13}}, {{-v11, -v12, -v13}, 0}};
dC1 = {{-v11, {u31, u32, u33}}, {{0, u13, -u12}, v11}};
dC2 = {{-v12, {u41, u42, u43}}, {{-u13, 0, u11}, v12}};
dC3 = {{-v13, {u51, u52, -u31 - u42}}, {{u12, -u11, 0}, v13}};
dD1 = {{-u11, {0, v13, -v12}}, {{-u31, -u41, -u51}, u11}};
dD2 = {{-u12, {-v13, 0, v11}}, {{-u32, -u42, -u52}, u12}};
dD3 = {{a8, {u81, u82, u83}}, {{v81, v82, v83}, b8}};
```

At last, the relation $C_1\ C_2 = D_3$

```
D3 == prod[C1, C2]
```

True

implies that

```
-dD3 + prod[dC1, C2] + prod[C1, dC2]
```

$$\begin{pmatrix} -a_8 - u_{13} & \{v_{12} - u_{81}, -u_{82} - v_{11}, -u_{83}\} \\ \{-u_{33} - v_{81}, -u_{43} - v_{82}, u_{31} + u_{42} - v_{83}\} & u_{13} - b_8 \end{pmatrix}$$

must be zero. Then we obtain $a_8 = -u_{13}$, $b_8 = u_{13}$, $u_{81} = v_{12}$, $u_{82} = -v_{11}$, $u_{83} = 0$, $v_{81} = -u_{33}$, $v_{82} = -u_{43}$ and $v_{83} = u_{31} + u_{42}$. We write then

```
dA = {{0, {u11, u12, u13}}, {{v11, v12, v13}, 0}};
dB = {{0, {-u11, -u12, -u13}}, {{-v11, -v12, -v13}, 0}};
dC1 = {{-v11, {u31, u32, u33}}, {{0, u13, -u12}, v11}};
dC2 = {{-v12, {u41, u42, u43}}, {{-u13, 0, u11}, v12}};
dC3 = {{-v13, {u51, u52, -u31 - u42}}, {{u12, -u11, 0}, v13}};
dD1 = {{-u11, {0, v13, -v12}}, {{-u31, -u41, -u51}, u11}};
dD2 = {{-u12, {-v13, 0, v11}}, {{-u32, -u42, -u52}, u12}};
dD3 = {{-u13, {v12, -v11, 0}}, {{-u33, -u43, u31 + u42}, u13}};
```

We finally obtain the 14 parameters and we can write

```
DerOs = {{0, 0, u11, u12, u13, v11, v12, v13}, {0, 0, -u11, -u12, -u13, -v11, -v12, -v13},
        {-v11, v11, u31, u32, u33, 0, u13, -u12}, {-v12, v12, u41, u42, u43, -u13, 0, u11},
        {-v13, v13, u51, u52, -u31 - u42, u12, -u11, 0},
  {-u11, u11, 0, v13, -v12, -u31, -u41, -u51},
        {-u12, u12, -v13, 0, v11, -u32, -u42, -u52},
  {-u13, u13, v12, -v11, 0, -u33, -u43, u31 + u42}}
```

$$
\begin{pmatrix}
0 & 0 & u_{11} & u_{12} & u_{13} & v_{11} & v_{12} & v_{13} \\
0 & 0 & -u_{11} & -u_{12} & -u_{13} & -v_{11} & -v_{12} & -v_{13} \\
-v_{11} & v_{11} & u_{31} & u_{32} & u_{33} & 0 & u_{13} & -u_{12} \\
-v_{12} & v_{12} & u_{41} & u_{42} & u_{43} & -u_{13} & 0 & u_{11} \\
-v_{13} & v_{13} & u_{51} & u_{52} & -u_{31}-u_{42} & u_{12} & -u_{11} & 0 \\
-u_{11} & u_{11} & 0 & v_{13} & -v_{12} & -u_{31} & -u_{41} & -u_{51} \\
-u_{12} & u_{12} & -v_{13} & 0 & v_{11} & -u_{32} & -u_{42} & -u_{52} \\
-u_{13} & u_{13} & v_{12} & -v_{11} & 0 & -u_{33} & -u_{43} & u_{31}+u_{42}
\end{pmatrix}
$$

with

```
Length[Variables[DerOs]]
```

14

We define then two new commands in order to prove that any matrix as the previous one is a derivation:

$$
\texttt{Coor}\left[\begin{pmatrix} \alpha_ & \texttt{x_} \\ \texttt{y_} & \beta_ \end{pmatrix}\right] := \{\alpha, \beta, \texttt{x[[1]]}, \texttt{x[[2]]}, \texttt{x[[3]]}, \texttt{y[[1]]}, \texttt{y[[2]]}, \texttt{y[[3]]}\}
$$

$$
\texttt{Zorn[x_]} := \begin{pmatrix} \texttt{x[[1]]} & \{\texttt{x[[3]]}, \texttt{x[[4]]}, \texttt{x[[5]]}\} \\ \{\texttt{x[[6]]}, \texttt{x[[7]]}, \texttt{x[[8]]}\} & \texttt{x[[2]]} \end{pmatrix}
$$

We must now check that $D(x\,y) = D(x)\,y + x\,D(y)$ for all $x, y \in \mathbf{O}_s$ . We write X and Y two generic elements in $\mathbf{O}_s$

```
X = α1 * A + α2 * B + Sum[λi * Ci, {i, 3}] + Sum[λi * Di-3, {i, 4, 6}]
```

$$
\begin{pmatrix} \alpha_1 & \{\lambda_1, \lambda_2, \lambda_3\} \\ \{\lambda_4, \lambda_5, \lambda_6\} & \alpha_2 \end{pmatrix}
$$

```
Y = β1 * A + β2 * B + Sum[μi * Ci, {i, 3}] + Sum[μi * Di-3, {i, 4, 6}]
```

$$
\begin{pmatrix} \beta_1 & \{\mu_1, \mu_2, \mu_3\} \\ \{\mu_4, \mu_5, \mu_6\} & \beta_2 \end{pmatrix}
$$

Then, we only have to compute

```
Simplify[Zorn[Coor[prod[X, Y]].DerOs] -
  prod[Zorn[Coor[X].DerOs], Y] - prod[X, Zorn[Coor[Y].DerOs]]]
```

$$
\begin{pmatrix} 0 & \{0, 0, 0\} \\ \{0, 0, 0\} & 0 \end{pmatrix}
$$

to prove our objective.

## ■ Basis of Der $(\mathbf{O}_s)$

It is easy now to compute a basis of the Lie algebra:

```
var = Variables[DerOs]; num = Length[var];
reglas = Table[var[[j]] -> δi,j, {i, num}, {j, num}];
Do[xi = DerOs //. reglas[[i]], {i, num}]
```

And then a generic element is as

$$\sum_{i=1}^{14} \alpha_i * x_i$$

$$\begin{pmatrix} 0 & 0 & \alpha_1 & \alpha_2 & \alpha_3 & \alpha_{12} & \alpha_{13} & \alpha_{14} \\ 0 & 0 & -\alpha_1 & -\alpha_2 & -\alpha_3 & -\alpha_{12} & -\alpha_{13} & -\alpha_{14} \\ -\alpha_{12} & \alpha_{12} & \alpha_4 & \alpha_5 & \alpha_6 & 0 & \alpha_3 & -\alpha_2 \\ -\alpha_{13} & \alpha_{13} & \alpha_7 & \alpha_8 & \alpha_9 & -\alpha_3 & 0 & \alpha_1 \\ -\alpha_{14} & \alpha_{14} & \alpha_{10} & \alpha_{11} & -\alpha_4 - \alpha_8 & \alpha_2 & -\alpha_1 & 0 \\ -\alpha_1 & \alpha_1 & 0 & \alpha_{14} & -\alpha_{13} & -\alpha_4 & -\alpha_7 & -\alpha_{10} \\ -\alpha_2 & \alpha_2 & -\alpha_{14} & 0 & \alpha_{12} & -\alpha_5 & -\alpha_8 & -\alpha_{11} \\ -\alpha_3 & \alpha_3 & \alpha_{13} & -\alpha_{12} & 0 & -\alpha_6 & -\alpha_9 & \alpha_4 + \alpha_8 \end{pmatrix}$$

We can now use **Recon** to compute the multiplication table of the Lie algebra

```
tabla = Array[a, {14, 14}];
Do[a[i, j] = Recon[c[xi , xj]], {i, 14}, {j, 14}];
tabla
```

$$\begin{pmatrix} 0 & -2x_{14} & 2x_{13} & x_1 & x_2 & x_3 & 0 & 0 & 0 & 0 & 0 & 2x_4 - x_8 & 3x_7 & 3x_{10} \\ 2x_{14} & 0 & -2x_{12} & 0 & 0 & 0 & x_1 & x_2 & x_3 & 0 & 0 & 3x_5 & 2x_8 - x_4 & 3x_{11} \\ -2x_{13} & 2x_{12} & 0 & -x_3 & 0 & 0 & 0 & -x_3 & 0 & x_1 & x_2 & 3x_6 & 3x_9 & -x_4 - x_8 \\ -x_1 & 0 & x_3 & 0 & x_5 & 2x_6 & -x_7 & 0 & x_9 & -2x_{10} & -x_{11} & x_{12} & 0 & -x_{14} \\ -x_2 & 0 & 0 & -x_5 & 0 & 0 & x_4 - x_8 & x_5 & x_6 & -x_{11} & 0 & 0 & x_{12} & 0 \\ -x_3 & 0 & 0 & -2x_6 & 0 & 0 & -x_9 & -x_6 & 0 & x_4 & x_5 & 0 & 0 & x_{12} \\ 0 & -x_1 & 0 & x_7 & x_8 - x_4 & x_9 & 0 & -x_7 & 0 & 0 & -x_{10} & x_{13} & 0 & 0 \\ 0 & -x_2 & x_3 & 0 & -x_5 & x_6 & x_7 & 0 & 2x_9 & -x_{10} & -2x_{11} & 0 & x_{13} & -x_{14} \\ 0 & -x_3 & 0 & -x_9 & -x_6 & 0 & 0 & -2x_9 & 0 & x_7 & x_8 & 0 & 0 & x_{13} \\ 0 & 0 & -x_1 & 2x_{10} & x_{11} & -x_4 & 0 & x_{10} & -x_7 & 0 & 0 & x_{14} & 0 & 0 \\ 0 & 0 & -x_2 & x_{11} & 0 & -x_5 & x_{10} & 2x_{11} & -x_8 & 0 & 0 & 0 & x_{14} & 0 \\ x_8 - 2x_4 & -3x_5 & -3x_6 & -x_{12} & 0 & 0 & -x_{13} & 0 & 0 & -x_{14} & 0 & 0 & -2x_3 & 2x_2 \\ -3x_7 & x_4 - 2x_8 & -3x_9 & 0 & -x_{12} & 0 & 0 & -x_{13} & 0 & 0 & -x_{14} & 2x_3 & 0 & -2x_1 \\ -3x_{10} & -3x_{11} & x_4 + x_8 & x_{14} & 0 & -x_{12} & 0 & x_{14} & -x_{13} & 0 & 0 & -2x_2 & 2x_1 & 0 \end{pmatrix}$$